\documentclass[12pt]{article}
\usepackage[intlimits]{amsmath}
\usepackage{amsfonts,amssymb,amscd,amsthm}
\input xy
\input epsf
\xyoption{all}

\newtheorem{theorem}{Theorem}[section]
\newtheorem{lemma}[theorem]{Lemma}
\newtheorem{proposition}[theorem]{Proposition}
\newtheorem{corollary}[theorem]{Corollary}
\theoremstyle{definition}

\title{Atiyah--Bott index on stratified manifolds\thanks{Research
supported in part by RFBR grants Nos.~05-01-00982 and~06-01-00098
and DFG grant 436 RUS 113/849/0-1\circledR ``$K$-theory and
noncommutative geometry of stratified manifolds." }}
\author{V.E.~Nazaikinskii, A.Yu.~Savin, and B.Yu.~Sternin}
\date{}
\begin{document}

\maketitle

\begin{abstract}
We define Atiyah--Bott index on stratified manifolds and express it
in topological terms. By way of example, we compute this index for
geometric operators on manifolds with edges.
\end{abstract}

\section{Introduction}

The paper deals with elliptic theory on stratified manifolds. The
symbol of a pseudodifferential operator on a stratified manifold is
a collection of symbols on the strata. The symbol on the stratum of
maximal dimension, called the \emph{interior symbol}, is of
particular importance, since it is a scalar function on the
cotangent bundle, while the symbols on the lower-dimensional strata
are operator-valued.

An operator with elliptic interior symbol defines an element in the
$K$-group of operators with zero interior symbol.\footnote{Similar
invariants were studied in other situations, e.g.,
in~\cite{Upm1,Mon1,MoNi1}.} This element, which we call the
\emph{Atiyah--Bott index}, has the following properties.
\begin{itemize}
    \item It is determined by
          the interior symbol of the operator
          and is a homotopy invariant of the interior symbol.
    \item For a smooth manifold, this invariant
          coincides with the Fredholm index.
    \item The Atiyah--Bott index is the
obstruction to making an operator with invertible interior symbol
invertible by adding lower-order terms.
    \item For a manifold with
boundary, this index coincides with the Atiyah--Bott obstruction
(see~\cite{AtBo2}) to the existence of well-posed (Fredholm)
boundary conditions for an elliptic operator in a bounded domain.
\end{itemize}

The main result of this paper is a general formula expressing the
Atiyah--Bott index in topological terms. We also compute the range
of the index mapping, i.e., the $K$-group of operators with zero
interior symbol.

This research was carried out during our stay at the Institute for
Analysis, Hannover University (Germany).  We are grateful to
Professor E.~Schrohe and other members of the university staff for
their kind hospitality.

\section{Atiyah--Bott index}

First, we recall the key properties of pseudodifferential operators
on stratified manifolds. For detailed exposition, e.g.,
see~\cite{NaSaSt4,NaSaSt5,NaSaSt3}.

\paragraph{Stratified manifolds.}

Let $\mathcal{M}$ be a compact stratified manifold in the sense
of~\cite{NaSaSt3}. Recall that $\mathcal{M}$ has a decreasing
filtration
\begin{equation*}
\label{filtr1} \mathcal{M}=\mathcal{M}_0\supset
\mathcal{M}_1\supset \mathcal{M}_2\ldots \supset
\mathcal{M}_N\supset \emptyset
\end{equation*}
of length $N$ by closed subsets  $\mathcal{M}_j$ such that the
complement $\mathcal{M}_{j}\setminus \mathcal{M}_{j+1}$ (an
\emph{open stratum}) is homeomorphic to the interior $M_j^\circ$ of
a compact manifold $M_j$ with corners (the \emph{blowup} of
$\mathcal{M}_{j}$). Let us denote the blowup of $\mathcal{M}$ by
$M$. In addition, any $x\in \mathcal{M}_{j}\setminus
\mathcal{M}_{j+1}$ has a neighborhood homeomorphic to the product
$$
U_x\times K_{\Omega_j},
$$
where $U_x$ is a neighborhood of $x$ in ${M}^\circ_{j}$ and
$$
K_{\Omega_j}=[0,1)\times\Omega_j\bigr/\{0\}\times \Omega_j
$$
is a cone whose base $\Omega_j$ is a stratified manifold with
filtration of length $<N$. In particular, $\mathcal{M}_N$ is a
smooth manifold.

\paragraph{Pseudodifferential operators on stratified manifolds.}

Let $ \Psi(\mathcal{M})$ be the algebra of pseudodifferential
operators of order zero on $\mathcal{M}$ acting in the space
$L^2(\mathcal{M})$ of complex-valued functions. (For the definition
of the algebra and of the measure defining the $L^2$-space, we
refer the reader to \cite{PlSe6} or \cite{NaSaSt4,NaSaSt5}.) The
Calkin algebra (the \emph{algebra of symbols}) is denoted by
$\Sigma(\mathcal{M})$.

The symbol $\sigma(D)$ of an operator $D$ on $\mathcal{M}$ is a
collection
$$
\sigma(D)=(\sigma_0(D),\sigma_1(D),...,\sigma_N(D))
$$
of symbols on the strata, where the symbol $\sigma_j(D)$ is defined
on the cosphere bundle $S^*M_j$ of the blowup of the corresponding
stratum. The symbol $\sigma_0(D)$ on the stratum
$\mathcal{M}\setminus \mathcal{M}_1$ of maximal dimension is a
scalar function called the \emph{interior symbol}, and the
remaining components of the symbol are operator-valued functions.

\paragraph{Definition of Atiyah--Bott index.}

Let
$$
\sigma_0\colon \Psi(\mathcal{M})\longrightarrow C(S^*M)
$$
be the mapping taking each operator to its interior symbol. This
mapping is surjective. Consider the short exact sequence
\begin{equation}
\label{short1} 0\longrightarrow J\longrightarrow
\Psi(\mathcal{M})\stackrel{\sigma_0}\longrightarrow
C(S^*M)\longrightarrow 0
\end{equation}
of $C^*$-algebras, where $J\subset\Psi(\mathcal{M})$ is the ideal
of operators with zero interior symbol.

The boundary mapping
\begin{equation} \label{bou1}
\delta\colon K_*(C(S^*M))\longrightarrow K_{*+1}(J)
\end{equation}
induced in $K$-theory by the exact sequence \eqref{short1} is
called the \emph{Atiyah--Bott index.}

Let us explain why the mapping \eqref{bou1} is called an index. If
$\mathcal{M}$ is a closed smooth manifold, then $M=\mathcal{M}$.
(The blowup coincides with the original manifold.) Moreover, the
kernel of $\sigma_0$ coincides with the ideal $\mathcal{K}$ of
compact operators, and the nontrivial boundary mapping
$$
\delta\colon K_1(C(S^*M))\longrightarrow K_{0}(\mathcal{K})\simeq
\mathbb{Z}
$$
is given by the Fredholm index.\footnote{Indeed, an element of
$K_1(C(S^*M))$ is determined by an invertible matrix with entries
in $C(S^*M)$. Consider the matrix as the symbol of some operator;
then the mapping takes (the equivalence class of) the matrix to the
Fredholm index of that operator.}

It turns out that the Atiyah--Bott index is the obstruction to
making an operator with invertible interior symbol invertible by
adding operators with zero interior symbol (cf.~\cite{MoNi1}).
\begin{theorem}\label{obst1}
Let $A$ be a matrix pseudodifferential operator with invertible
interior symbol $\sigma_0(A)$ on a stratified manifold
$\mathcal{M}$. A necessary and sufficient condition that there
exists an operator $R$ with zero interior symbol such that $A+R$ is
invertible is that the Atiyah--Bott index $\delta[\sigma_0(A)]\in
K_0(J)$ is zero.
\end{theorem}

The proof is given in the Appendix.

\section{Main theorem}

\paragraph{Analytic $K$-homology.}

Let us recall several facts about analytic  $K$-homology. (Detailed
exposition and further references can be found
in~\cite{HiRo1},~Chapter~5.)

Let $\mathcal{M}$ be a compact stratified manifold. By
$\mathcal{D}(\mathcal{M})$ we denote the algebra of  \emph{local
operators} on $L^2(\mathcal{M})$, i.e., operators compactly
commuting with multiplications by continuous functions on
$\mathcal{M}$. If $\mathcal{M}'$ is a closed subspace, then by
$\mathcal{D}(\mathcal{M},\mathcal{M}')\subset
\mathcal{D}(\mathcal{M})$ we denote the ideal of \emph{locally
compact operators}. By definition, these are operators whose
composition with multiplications by continuous functions vanishing
on $\mathcal{M}'$ is compact.

The $K$-homology groups of $\mathcal{M}$, $\mathcal{M}'$, and
$\mathcal{M}\setminus\mathcal{M}'$ can be defined as
\begin{equation}\label{eqq1}
    K_*(\mathcal{M})\simeq
K_{*+1}\Bigl(\mathcal{D}(\mathcal{M})/\mathcal{K}\Bigr),\quad
K_*(\mathcal{M}')\simeq
K_{*+1}\Bigl(\mathcal{D}(\mathcal{M},\mathcal{M}')/\mathcal{K}\Bigr),
\end{equation}
\begin{equation}\label{eqq2}
K_*(\mathcal{M}\setminus\mathcal{M}')\simeq
K_{*+1}\Bigl(\mathcal{D}(\mathcal{M})/\mathcal{D}(\mathcal{M},\mathcal{M}')\Bigr).
\end{equation}
Moreover the  $K$-homology sequence for the pair
$\mathcal{M}'\subset \mathcal{M}$
$$
\ldots \to K_*(\mathcal{M}')\to K_*(\mathcal{M})\to
K_*(\mathcal{M}\setminus \mathcal{M}')\stackrel\partial\to
K_{*+1}(\mathcal{M}')\to\ldots
$$
coincides with the $K$-theory sequence for the short exact sequence
$$
0\to
\mathcal{D}(\mathcal{M},\mathcal{M}')/\mathcal{K}\longrightarrow
\mathcal{D}(\mathcal{M})/\mathcal{K}\longrightarrow\mathcal{D}(\mathcal{M})\bigl/
\mathcal{D}(\mathcal{M},\mathcal{M}')\to 0
$$
of $C^*$-algebras.

\paragraph{Main theorem.}

Pseudodifferential operators on stratified manifolds are local
(see~\cite{NaSaSt5}). On the other hand, operators in  $J$ are
locally compact with respect to the subspace
$\mathcal{M}_1\subset\mathcal{M}$; i.e., their compositions with
functions vanishing on $\mathcal{M}_1$ are compact. Thus, we have
the commutative diagram
\begin{equation}
\label{diag1}
\begin{array}{rcccccl}
  0\to & J & \to  & \Psi(\mathcal{M}) & \to & C(S^*M) & \to 0 \vspace{2mm}\\
   & i^{\mathcal{M}_1}\downarrow\phantom{i^{\mathcal{M}_1}} &  & \downarrow &  & \downarrow i^{M^\circ} &  \vspace{2mm}\\
  0\to & \mathcal{D}(\mathcal{M},\mathcal{M}_1) & \to & \mathcal{D}(\mathcal{M}) &
  \to & \mathcal{D}(\mathcal{M})\bigr/ \mathcal{D}(\mathcal{M},\mathcal{M}_1) & \to 0. \\
\end{array}
\end{equation}
In view of the isomorphisms \eqref{eqq1} and \eqref{eqq2}, the
vertical mappings in \eqref{diag1} induce homomorphisms of
$K$-groups into $K$-homology groups:
\begin{align*}
i^{\mathcal{M}_1}_*&\colon K_{*+1}(J)\longrightarrow
K_{*+1}\Bigl(\mathcal{D}(\mathcal{M},\mathcal{M}_1)/\mathcal{K}\Bigr)\simeq
K_*(\mathcal{M}_1),
\\
i^{M^\circ}_*&\colon K_{*}(C(S^*M))\longrightarrow
K_{*}\Bigl(\mathcal{D}(\mathcal{M})/\mathcal{D}(\mathcal{M},\mathcal{M}_1)\Bigr)\simeq
K_{*+1}(\mathcal{M}\setminus \mathcal{M}_1).
\end{align*}
\begin{theorem}
\label{mainth1} If $\mathcal{M}$ has no closed smooth components,
then the diagram
\begin{equation*}
\label{abformula} {{\xymatrix{
  K_*(C(S^*M))\ar[d]_{i^{M^\circ}_*}\ar[r]^\delta
 & K_{*+1}(J)\ar[d]^{i^{\mathcal{M}_1}_*}\\
 K_{*+1}(\mathcal{M}\setminus
\mathcal{M}_1)\ar[r]^\partial & \widetilde{K}_{*}(\mathcal{M}_1)
}}}
\end{equation*}
commutes. Moreover,  $i^{\mathcal{M}_1}_*$ is an isomorphism. Here
 \begin{itemize}
    \item $\widetilde{K}_{*}(\mathcal{M}_1)\simeq
     \ker\{K_*(\mathcal{M}_1)\to K_*(pt)\}$ is
     the reduced $K$-homology group.
    \item $\partial$ is the boundary mapping in exact $K$-homology
    sequence of the pair    $\mathcal{M}_1\subset \mathcal{M}.$
 \end{itemize}
\end{theorem}

\begin{proof}

1. Consider the commutative diagram
$$
\xymatrix{ 0 \ar[r] & J \ar[d]_j \ar[r] & \Psi \ar[d] \ar[r] &
C(S^*M)\ar@{=}[d]
\ar[r] & 0\\
0 \ar[r] & J/\mathcal{K} \ar[r] & \Sigma \ar[r] & C(S^*M) \ar[r] &
0,
}
$$
where $\Psi$ is the algebra of pseudodifferential operators,
$\Sigma=\Psi/\mathcal{K}$ is the algebra of symbols, and $j$ is the
natural projection. The boundary mappings corresponding to the
upper and lower rows of the diagram are compatible. Since $j$
induces a monomorphism in $K$-theory (see the lemma below), it
suffices to compute the boundary mapping corresponding to the lower
exact sequence.

\begin{lemma}\label{leka1}
The natural projection $j\colon J\to J/\mathcal{K}$ induces
isomorphisms
\begin{equation*}
\label{kth1} K_0(J)\simeq K_0(J/\mathcal{K}),\qquad K_1(J)=\ker
\bigl(\operatorname{ind}\colon K_1(J/\mathcal{K})\longrightarrow
\mathbb{Z}\bigr).
\end{equation*}
\end{lemma}
\begin{proof}
The boundary mapping in the $K$-theory exact sequence of the pair
$J\to J/\mathcal{K}$ is the index mapping
$$
K_1(J/\mathcal{K}) \longrightarrow K_0(\mathcal{K})=\mathbb{Z}.
$$
It is surjective. The exact sequence of pair $\mathcal{K}\subset J$
gives the desired isomorphism.
\end{proof}

2. Consider the commutative diagram
\begin{equation}\label{comm}
    \begin{array}{rcccccl}
      0\to& J/\mathcal{K} & \to &\Sigma&\to & C(S^*M)&\to 0 \\
      &\downarrow & &\downarrow & & \downarrow& \\
      0\to &\mathcal{D}(\mathcal{M}, \mathcal{M}_1)/\mathcal{K} &
      \to&
\mathcal{D}(\mathcal{M})/\mathcal{K} &\to  &
  \mathcal{D}(\mathcal{M})\bigr/ \mathcal{D}(\mathcal{M}, \mathcal{M}_1)&\to 0.\\
    \end{array}
\end{equation}
As we mentioned earlier, the sequence of $K$-groups induced by the
lower row in~\eqref{comm} is isomorphic to the $K$-homology exact
sequence of the pair $\mathcal{M}_1\!\subset \mathcal{M}$.

Since $\mathcal{M}$ is assumed to have no closed smooth components,
there exists a nonsingular vector field (a section $M\to S^*M$ of
the cosphere bundle) on the blowup $M$. It follows that
$K_*(\Sigma)$ and $K_*(C(S^*M))$ have direct summands $K_*(C(M))$,
which are mapped isomorphically onto each other by
$\sigma_0\colon\Sigma\to C(S^*M)$. Hence \eqref{comm} gives the
diagram
\begin{equation}\label{maindiag1}
 \begin{array}{ccccccc}
    \!\!\!K_*\!(J/\mathcal{K}) \to & \!\!\!\!\!K_*\!(\Sigma)/K_*\!(C(\!M))
   \!\!\!\!&\to  &
   \!\!\!\!\!K_*\!(C(S^*\!M))/K_*\!(C(\!M))\!\!\!\! &\!\!\to\! &\!\!\! K_{*\!+\!1}\!(J/\mathcal{K}) &
   \!\!\!\!\! \to \vspace{1mm}\\
      i_*^{\mathcal{M}_1}\downarrow\phantom{i_*^{\mathcal{M}_1}}   & i_*^\mathcal{M}\downarrow  & &
     i^{M^\circ}_*\downarrow & & i_*^{\mathcal{M}_1}\downarrow &
     \vspace{1mm} \\
 \!\!K_{*\!+\! 1}(\mathcal{M}_1\!)\!\!     \to & K_{* +
1}(\mathcal{M}) &\to & K_{*+1}(\mathcal{M}\setminus \mathcal{M}_1
)&\!\stackrel\partial \to\!& K_{*}(\mathcal{M}_1) &
    \!\!\!\!\!\to\\
 \end{array}
\end{equation}
of $K$-groups with exact rows.

We proved in  \cite{NaSaSt3} that $i_*^\mathcal{M}$ and
$i^{M^\circ}_*$ are isomorphisms. Hence $i_*^{\mathcal{M}_1}$ is an
isomorphism as well by the five lemma. Now the commutative
rightmost square in \eqref{maindiag1}, together with
Lemma~\ref{leka1}, gives the desired commutative diagram
$$
\xymatrix{
 K_*(C(S^*M))\ar[d]_{i_*^{M^\circ}}\ar[r]^\delta
 & K_{*+1}(J)\ar[d]^{i_*^{\mathcal{M}_1}}\\
 K_{*+1}(\mathcal{M}\setminus \mathcal{M}_1)\ar[r]^\partial &
 \widetilde{K}_{*}(\mathcal{M}_1)
}
$$
where $i_*^{\mathcal{M}_1}$ is an isomorphism.

The proof of Theorem~\ref{mainth1} is complete.
\end{proof}

\paragraph{Cohomological index formula.}

To be definite, we write out a (co)homological formula for an
element $[a]\in K_1(C(S^*M))$ represented by an invertible symbol
$a$ on $S^*M$. Consider the $K$-homology element
\begin{equation}\label{qq3}
    i_*^{\mathcal{M}_1}\delta[a]\in K_1(\mathcal{M}_1).
\end{equation}
The Chern character of this element is a rational homology class on
$\mathcal{M}_1$. Let us evaluate the pairing of this class with an
arbitrary cohomology class on $\mathcal{M}_1$.

To this end, let $Y\subset M$ be a compact smooth manifold  of
dimension ${\rm dim} M-1 $,  homeomorphic to  the boundary
$\partial M=M\setminus M^\circ$. (The smooth structure is obtained
by smoothing the corners.) Next, let
$$
\pi\colon  Y\to \mathcal{M}_1
$$
be the projection onto the singularity set $\mathcal{M}_1$.

\begin{corollary}\label{sled1}
For any $[a]\in K_{1}(C(S^*M))$ and $x\in H^{odd}(\mathcal{M}_1)$,
one has
\begin{equation}
\label{coho} \Bigl\langle\operatorname{ch}
\bigl(i_*^{\mathcal{M}_1}\delta[a]\bigr),x\Bigr\rangle=
{\Bigl\langle\operatorname{ch} \left[a|_Y\right]\cdot
\operatorname{Td}(T^*Y\otimes\mathbb{C})\pi^*x,[T^*Y\times
\mathbb{R}]\Bigr\rangle},
\end{equation}
where $\operatorname{ch} (i_*^{\mathcal{M}_1}\delta[a])\in
H_{odd}(\mathcal{M}_1,\mathbb{Q})$ is the homological Chern
character and $[a|_Y]\in K^0_c(T^*Y\times \mathbb{R})$ is the
restriction of interior symbol $a$ to $Y$.
\end{corollary}

\begin{proof}
Since the Chern character is a rational isomorphism of $K$-theory
and homology, it suffices to evaluate the pairing with an element
$x$ of the form $x=\operatorname{ch} y$, where $y\in
K^1(\mathcal{M}_1)$.

By Theorem~\ref{mainth1},
$$
\Bigl\langle\operatorname{ch}
\bigl(i_*^{\mathcal{M}_1}\delta[a]\bigr),\operatorname{ch}
y\Bigr\rangle= \Bigl\langle i_*^{\mathcal{M}_1}\delta[a],y
\Bigr\rangle= \Bigl\langle \partial i_*^{M^\circ}[a],
y\Bigr\rangle.
$$
It follows from the properties of the boundary mapping $\partial$
in $K$-homology that
$$
\Bigl\langle \partial \cdot i_*^{M^\circ}[a], y\Bigr\rangle=
\Bigl\langle \pi_*\cdot \partial_0 \cdot i_*^{M^\circ}[a],
y\Bigr\rangle,
$$
where $\partial_0\colon  K_0(\mathcal{M}\setminus \mathcal{M}_1)\to
K_1(Y)$ is the boundary mapping for the pair $\partial M\subset M$.
We now transfer the computation of the pairing to $Y$:
$$
\Bigl\langle\pi_* \partial_0  i_*^{M^\circ}[a], y\Bigr\rangle=
\Bigl\langle \partial_0  i_*^{M^\circ}[a], \pi^*y\Bigr\rangle.
$$
By applying the Atiyah--Singer formula on $Y$, we obtain
$$
\Bigl\langle  \partial_0  i_*^{M^\circ}[a], \pi^*y\Bigr\rangle=
\Bigl\langle\operatorname{ch} \left[a|_Y\right]\cdot
\operatorname{Td}(T^*Y\otimes\mathbb{C})\operatorname{ch}\pi^*y,[T^*Y\times
\mathbb{R}]\Bigr\rangle.
$$
The right-hand side coincides with the desired cohomological
expression~\eqref{coho}, since $x=\operatorname{ch} y$.
\end{proof}

\section{Examples}

In this section, we compute the Atiyah--Bott index for geometric operators on
stratified manifolds with stratification of length one. Such manifolds are
called \emph{manifolds with edges}. Geometrically, a manifold with edges is
obtained as follows. Take a smooth manifold $M$ whose boundary $\partial M$ is
fibered over a smooth base $X$, $\pi\colon \partial M\to X$. Then identify the
points in each fiber of $\pi$. What is obtained is a manifold $\mathcal{M}$
with edge $X$. The blowup of $\mathcal{M}$ is just $M$, and the manifold $Y$ in
Corollary~\ref{sled1} is diffeomorphic to the boundary $\partial M$.

The Atiyah--Bott index is zero for the Beltrami--Laplace and Euler
operators, since the restrictions of their principal symbols to $Y$
give rise to zero elements in $K^0_c(T^*Y\times\mathbb{R})$.

Let us (rationally) compute the Atiyah--Bott index of the Dirac
operator.

Suppose that $M$,  $X$, and $\pi$ are equipped with spin
structures. Next, assume that the induced spin structure on the
total space $\partial M$ of the bundle $\pi$ is compatible with the
spin structure on $\partial M$. Let
\begin{equation*}
\mathcal{D}\colon S_+(M^\circ)\longrightarrow S_-(M^\circ)
\end{equation*}
be the Dirac operator on $M^\circ$. (We assume that the dimension
of $M$ is even.)

\begin{proposition}
For the Dirac operator $\mathcal{D}$ on a manifold $\mathcal{M}$
with edge $X$, the homology class
$$
\operatorname{ch}(i_*^X\delta[\sigma_0(\mathcal{D})])\in
H_{*}(X)\otimes \mathbb{Q}
$$
is Poincar\'e dual to the cohomology class
$$
A(X)\cdot \pi_*(A(\Omega)),
$$
where $A(X)$ and  $A(\Omega)$ are the $A$-classes of the base and
of the fibers of $\pi$, respectively, and $\pi_*\colon H^*(Y)\to
H^{*-\dim\Omega}(X)$ stands for integration over the fibers.
\end{proposition}
\begin{proof}
Choosing a connection in $\pi$, we obtain a decomposition
\begin{equation*}
\label{dirsum1} T^* M|_Y\simeq \pi^* T^*X\oplus  (T^*\Omega\oplus
\mathbb{R})
\end{equation*}
of the restriction of the cotangent bundle to $Y$ into horizontal
and vertical components. With regard to this decomposition, the
symbol of the Dirac operator over $Y$ is the tensor product of the
pull-back of the symbol of the Dirac operator on the base by the
family of Dirac operators in the fibers. This gives a factorization
$$
[\sigma(\mathcal{D})|_Y]=\bigl(\pi^*[\sigma(\mathcal{D}_X)]\bigr)[\sigma(\mathcal{D}_{\Omega\times
\mathbb{R}})],
$$
of the corresponding element in $K$-theory, where
$[\sigma_0(\mathcal{D}_X)]\in K^{\dim X}_c(T^*X)$ and
$[\sigma(\mathcal{D}_{\Omega\times \mathbb{R}})]\in
K^{\dim\Omega+1}_c(T^*\Omega\oplus \mathbb{R})$.

Let $x\in H^{odd}(X)$ be a cohomology class. Let us compute the
pairing
$$
\bigl<\operatorname{ch}(i_*^X\delta[\sigma_0(\mathcal{D})]),x\bigr>.
$$
By Corollary~\ref{sled1}, this number is equal to
$$
\bigl<\operatorname{ch}[\sigma_0(\mathcal{D})|_Y]\operatorname{Td}(T^*Y\otimes
\mathbb{C})\pi^*x,[T^*Y\times \mathbb{R}]\bigr>.
$$
Using the factorization of $[\sigma_0(\mathcal{D})|_Y]$, we see
that this expression is equal to
\begin{multline*}
\bigl<\operatorname{ch}[\sigma(\mathcal{D}_X)]\operatorname{ch}
[\sigma(\mathcal{D}_{\Omega\times \mathbb{R}})]
\operatorname{Td}(T^*X\otimes
\mathbb{C})\operatorname{Td}(T^*\Omega\otimes
\mathbb{C})\pi^*x,[T^*Y\times \mathbb{R}]\bigr>\\
=\bigl<\operatorname{ch}[\sigma(\mathcal{D}_X)]\pi_*\bigl(\operatorname{ch}
[\sigma(\mathcal{D}_{\Omega\times
\mathbb{R}})]\operatorname{Td}(T^*\Omega\otimes \mathbb{C})\bigr)
\operatorname{Td}(T^*X\otimes \mathbb{C}) x,[T^*X]\bigr>
\\
=\bigl<A(X)\pi_*(A(\Omega)) x,[X]\bigr>.
\end{multline*}
(At the last step, we have used the standard transition from
cohomology classes on the cotangent bundle to cohomology classes on
the manifolds themselves by using the Thom isomorphism.)

Thus,
$$
\bigl<\operatorname{ch}(i_*^X\delta[\sigma_0(\mathcal{D})]),x\bigr>=\bigl<A(X)\pi_*(A(\Omega))
x,[X]\bigr>
$$
for all $x\in H^{odd}(X)$. It follows that the classes
$\operatorname{ch}(i_*^X\delta[\sigma_0(\mathcal{D})])$ and
$A(X)\pi_*(A(\Omega))$ are Poincar\'e dual.
\end{proof}

A similar computation of the Atiyah--Bott index can be carried out
for the signature operator. In this case, one should replace the
$A$-classes by the $L$-classes.

\section*{Appendix. Proof of Theorem~\ref{obst1}}

\paragraph{Auxiliary lemmas.}
The following two lemmas are standard.
\begin{lemma}\label{lem1}
Let $A\colon\mathcal{H}_1\to \mathcal{H}_2$ be a continuous family
of Fredholm operators acting in spaces of sections of
infinite-dimensional Hilbert bundles over a locally compact base
$X$. Suppose that  $A$ is invertible at infinity. For the existence
of a continuous family of finite rank operators $R(x)$, $x\in X$,
vanishing at infinity such that the family $A+R$ is everywhere
invertible, it is necessary and sufficient that
$$
\operatorname{ind} A=0\in K^0_c(X).
$$
\end{lemma}

\begin{lemma}\label{lem2}
Let $\mathcal{H}$ be an infinite-dimensional Hilbert bundle over a
locally compact base $X$. Then each element of
$K_1(C_0(X,\mathcal{K}(\mathcal{H})))$ has a representative that is
an invertible element of the unital algebra
$C_0(X,\mathcal{K}(\mathcal{H}))^+$.
\end{lemma}

\paragraph{Proof of Theorem~\ref{obst1}.}

If there exists an $R$ with the desired properties, then $A+R$ is
invertible, $[A+R]\in K_1(\Psi(\mathcal{M}))$, and
$(\sigma_0)_*[A+R]=[\sigma_0(A)]$. We conclude that
$\delta[\sigma_0(A)]=\delta((\sigma_0)_*[A+R])=0$, since the
sequence
$$
\ldots\to K_1(\Psi)\stackrel{(\sigma_0)_*}\to
K_1(C(S^*M))\stackrel\delta\to K_0(J)\to \ldots
$$
is exact.

Let us prove the converse. To this end, let us recall several
properties of pseudodifferential operators on stratified manifolds
(see \cite{NaSaSt4,NaSaSt5}).

The algebra $\Psi=\Psi(\mathcal{M})$ is solvable of length $(2N+1)$
(e.g., see \cite{NaSaSt4}), where $N$ is the length of the
stratification of $\mathcal{M}$, with composition series
\begin{equation}
\label{compo1}
\Psi=\Psi_0\supset\Psi_1\supset\Psi_2\supset\Psi_3\supset\ldots\supset
\Psi_{2N}\supset \Psi_{2N+1}=\mathcal{K}\supset \{0\}.
\end{equation}
To describe the ideals $\Psi_j$ in more detail, recall that the
symbol of an operator $D$ is the collection
$$
\sigma(D)=(\sigma_0(D),\sigma_1(D),...,\sigma_N(D))
$$
of symbols on the strata. Here $\sigma_j(D)$ is a family of
operators in $L^2(K_{\Omega_j})$, where $K_{\Omega_j}$ is cone with
base $\Omega_j$, parametrized by cosphere bundle  $S^*M_j$; the
conormal symbol (corresponding to the cone tips) of this family
does not depend on the covariables in $S^*M_j$ and is a family,
parametrized by $\mathbb{R}\times M_j$, of pseudodifferential
operators on the base $\Omega_j$ of the cone.

Now the ideals in~\eqref{compo1} can be described as follows. The
ideal
$$
\Psi_{2j+1}\subset \Psi_{2j},\quad j\ge 0,
$$
consists of operators with zero symbol $\sigma_j$ on the open
stratum $\mathcal{M}^\circ_j$, and the ideal
$$
\Psi_{2j}\subset \Psi_{2j-1},\quad j\ge 1,
$$
consists of operators with zero conormal symbol
$\sigma_c(\sigma_j)$ of $\sigma_j$ on $\mathcal{M}_j^\circ$.
Moreover, we have the isomorphisms
$$
\Psi_{2j}/\Psi_{2j+1}\simeq C(S^*M_j,\mathcal{K}L^2(K_{\Omega_j}))
$$
(the mapping is defined by the symbol $\sigma_j$) and
$$
\Psi_{2j-1}/\Psi_{2j}\simeq C_0(\mathbb{R}\times
M_j,\mathcal{K}L^2(\Omega_j))
$$
(the mapping is defined by the conormal symbol
$\sigma_c(\sigma_j)$).

To construct an invertible perturbation of an operator with trivial
Atiyah--Bott index, we use Proposition~\ref{prop1} below for
$j=1,2,\dotsc,2N+1$. As a result, we shall obtain an invertible
operator with interior symbol equal to that of the original
operator. This will end the proof of Theorem~\ref{obst1}.

Let $A\colon L^2(\mathcal{M},\mathbb{C}^n)\to
L^2(\mathcal{M},\mathbb{C}^n)$ be a matrix pseudodifferential
operator. We say that it is \emph{invertible modulo the ideal}
$\Psi_j$ if there exists a matrix operator $B$ such that the
compositions $AB$ and $BA$ are equal to identity modulo operators
with matrix entries in $\Psi_j$.

\begin{proposition}\label{prop1}
Let $A$ be a matrix pseudodifferential operator on $\mathcal{M}$
such that for some $j$ it is invertible modulo ideal $\Psi_j$ and
$$
\delta_j[A]=0,
$$
where $\delta_j\colon K_1(\Psi/\Psi_j)\longrightarrow K_0(\Psi_j)$
is the boundary map in $K$-theory for the pair $\Psi \to
\Psi/\Psi_j$. Then there exists a pseudodifferential operator
$\widetilde{A}$ on $\mathcal{M}$ that is invertible modulo the
ideal $\Psi_{j+1}$, is equal to $A$ modulo $\Psi_j$, and has zero
index $\delta_{j+1}[\widetilde{A}]=0$ in $K_0(\Psi_{j+1})$.
\end{proposition}

\begin{proof}
Let $s_{2j}=\sigma_j$ and $s_{2j-1}=\sigma_c(\sigma_j)$. (In this
notation, the ideal $\Psi_{k+1}$ in $\Psi_k$ is determined by the
condition $s_{k}=0$.)

0. Consider the diagram
$$
\xymatrix{ K_1(\Psi_{j}/\Psi_{j+1})\ar@{=}[d]\ar[r] &
K_1(\Psi/\Psi_{j+1})\ar[r] \ar[d]^{\delta_{j+1}} &
K_1(\Psi/\Psi_{j})\ar[r]^{\delta''}\ar[d]^{\delta_j} &
 K_0(\Psi_{j}/\Psi_{j+1})\ar@{=}[d]\\
K_1(\Psi_{j}/\Psi_{j+1})\ar[r]^{\delta'} &
K_0(\Psi_{j+1})\ar[r]^\gamma &
K_0(\Psi_{j})\ar[r] & K_0(\Psi_{j}/\Psi_{j+1})\\
}
$$
The middle columns are parts of the exact sequences for the pairs
$\Psi\to\Psi/\Psi_{j+1}$ and $\Psi\to \Psi/\Psi_j$. The rows are
the exact sequences of the pairs  $\Psi/\Psi_{j+1}\to
\Psi/\Psi_{j}$ and $\Psi_j\to \Psi_j/\Psi_{j+1}$. Since the
boundary mapping is natural, it follows that the diagram commutes.

1. Suppose that $\delta_j[A]=0$. Then $\delta''[A]=0$. Hence the
element
$$
\operatorname{ind} s_{j+1}=\delta''[A]\in K_0(\Psi_j/\Psi_{j+1})
$$
is zero. (Indeed, $s_{j+1}$ is a Fredholm family invertible at
infinity, since $s_j$ is invertible by assumption.) Thus, by
Lemma~\ref{lem1}, there exists a pseudodifferential operator
$B=A\mod\Psi_{j}$ on $\mathcal{M}$ invertible modulo $\Psi_{j+1}$.

2. Since $\gamma\delta_{j+1}[B]=\delta_j[A]=0$, we have
$\delta_{j+1}[B]=\delta'(z)$ for some $z\!\in\!
K_1(\Psi_j/\Psi_{j+1})$.

3. For $\widetilde{A}$ we take the composition
$$
\widetilde{A}=BZ^{-1},
$$
where $Z\in \Psi_{j}^+$ is a pseudodifferential operator invertible
modulo $\Psi_{j+1}$ such that $[Z]=z$. (The existence of $Z$ is
guaranteed by Lemma~\ref{lem2}.)

It is clear that $\widetilde{A}$ has the desired properties: it is
invertible modulo $\Psi_{j+1}$, is equal to $A$ modulo $\Psi_j$,
and has zero index $\delta_{j+1}[\widetilde{A}]=0$, since
$$
\delta_{j+1}[\widetilde{A}]=\delta_{j+1}[B]-\delta_{j+1}[Z]=
\delta_{j+1}[B]-\delta'(z)=0
$$
by construction.

The proof of  Proposition~\ref{prop1} is complete.

\end{proof}



\end{document}